\documentclass{article}

\usepackage{amsmath,amsthm}     
\usepackage{url}
\usepackage{amsfonts} 

\date{}

\theoremstyle{theorem}
\newtheorem{theorem}{Theorem}

\theoremstyle{definition}

\begin{document}

\title{Generalizations of Bertrand's Postulate \\
to Sums of Any Number of Primes}

\author{Joel E. Cohen\\             
The Rockefeller University\\    
1230 York Avenue, Box 20\\
New York, NY 10065, USA\\
\url{cohen@rockefeller.edu}\\
Additional affiliations: Columbia University, University of Chicago
\\ 
Accepted for publication in \emph{Mathematics Magazine}\\ MATHMAG-D-21-00038R1 on September 03, 2021}                 

\maketitle

\noindent  In 1845, Bertrand \cite{bertrand} conjectured what became known as Bertrand's postulate:
twice any prime strictly exceeds the next prime.
Tchebichef \cite{Chebyshev1852} (to use the spelling he used on the original publication)
presented his proof of Bertrand's postulate to the Imperial Academy of St. Petersburg in 1850
and published it in 1852.
It is now sometimes called the Bertrand-Chebyshev theorem.

Surprisingly, a stronger statement seems not to be well known, 
but is elementary to prove:
the sum of any two
consecutive primes strictly exceeds the next prime,
except for the only equality $2+3=5$. 
After I conjectured and proved this statement independently, 
a very helpful referee pointed out 
that Ishikawa \cite{Ishikawa1934} published this result in 1934 
(with a different proof).

This observation is a special case of a much more general result, Theorem \ref{thm:general},
that is 
also elementary to prove (given the prime number theorem), and perhaps not previously noticed:
if
$p_n$ denotes the $n$th prime, $n=1, 2, \ldots$, with $p_1=2, p_2=3, \ldots$, and
if $c_1, \ldots, c_g$ are nonnegative integers (not necessarily distinct), and
$d_1, \ldots, d_h$ are positive integers (not necessarily distinct),
and $g>h\ge 1$, then there exists a positive integer $N$ 
such that $p_{n-c_1}+p_{n-c_2}+\cdots +p_{n-c_g}>p_{n+d_1}+\cdots +p_{n+d_h}$
for all $n\ge N$.
We prove this result using only the prime number theorem.
For any instance of this result,
we sketch a way to find the least possible $N$.

Bertrand's postulate $p_n+p_n\ge p_{n+1}$ is the special case
of Theorem \ref{thm:general} in which
$g=2, c_1=c_2=0, h=1, d_1=1,N=1$.

We give an easy, independent proof of some other special cases, 
notably Theorem \ref{thm:c1}:
for all $n>1, p_{n-1} + p_{n}\ge p_{n+1}$ and 
equality holds only for $n=2$. 
We give some numerical results and unanswered questions.

\section{Main result}

\begin{theorem}\label{thm:general} 
If $c_1, \ldots, c_g$ are $g>1$ nonnegative integers (not necessarily distinct),
and
$d_1, \ldots, d_h$ are $h$ positive integers (not necessarily distinct),
with $1 \leq h < g$,
then there exists a positive integer $N$ 
such that, for all $n\ge N$,
$$p_{n-c_1}+p_{n-c_2}+\cdots +p_{n-c_g}>p_{n+d_1}+\cdots +p_{n+d_h}.$$
\end{theorem}
\begin{proof}
For real-valued functions $f,\phi$ with real arguments $x$
such that $\phi(x)>0$ for all sufficiently large $x$, 
we define $f(x)\sim \phi(x)$ to mean that $\lim_{x\to\infty} f(x)/\phi(x)=1$.
The prime number theorem says that if $x>0$ and $\pi(x)$ is the number of primes that do not exceed $x$, then
$$\pi(x)\sim \frac{x}{\log x}.$$
Exactly $n$ primes do not exceed the $n$th prime $p_n$, 
so as $n\to\infty$,
\begin{align*}
\pi(p_n)=n&\sim \frac{p_n}{\log p_n},\\
\log n &\sim  \log{p_n}-\log\log p_n,\\
n\log n &\sim \frac{p_n}{\log p_n} \log{p_n}-\frac{p_n}{\log p_n}\log\log p_n
=p_n\left(1-\frac{\log\log p_n}{\log p_n}\right)\sim p_n.
\end{align*}
In summary, $p_n\sim n \log n$.
Therefore, for any fixed integer $C$ such that $p_{n\pm C}$ is defined, 
$p_{n\pm C}\sim (n\pm C)\log (n\pm C)\sim n\log n$. So
\begin{align*}
p_{n-c_1}+p_{n-c_2}+\cdots +p_{n-c_g}&\sim gn\log n,\\
p_{n+d_1}+\cdots +p_{n+d_h}&\sim hn\log n.
\end{align*}
Hence
\begin{equation*}
\frac{p_{n-c_1}+p_{n-c_2}+\cdots +p_{n-c_g}}{p_{n+d_1}+\cdots +p_{n+d_h}}\sim \frac{g}{h}>1.
\end{equation*}
Therefore there exists a positive integer $N$ 
such that 
$p_{n-c_1}+p_{n-c_2}+\cdots +p_{n-c_g}>p_{n+d_1}+\cdots +p_{n+d_h}$
for all $n\ge N$.
\end{proof}

\section{Finding the least $N$: an alternative approach}

The same helpful referee suggested that,
for greater specificity about  $N$ in Theorem \ref{thm:general},
these inequalities could also be proved using  
inequalities of Rosser and Schoenfeld \cite[p. 69, (3.12), (3.13)]{RosserSchoenfeld1962}:
\begin{align*}
n \log n < &p_n \quad  &for~ 1 \leq n,\\
&p_n < n(\log n + \log \log n)  \quad &for~ 6 \leq n.
\end{align*}
We sketch the idea.
Suppose we want to find the least positive integer $N$ 
such that
$$p_{n}+p_{n-1}+p_{n-2}>p_{n+1}+p_{n+2}
 \quad \quad \quad \quad for~ n\geq N.$$
We look numerically for the least $n\geq 8$  such that
\begin{align*}
&n \log n + (n-1) \log(n-1) + (n-2) \log(n-2) \\
&>(n+1)(\log(n+1) + \log \log(n+1)) + (n+2)(\log(n+2) + \log \log(n+2)).
\end{align*}
Then $n-2\geq 6$
and the Rosser-Schoenfeld inequalities apply to all terms above.
Such an $n$ must exist because 
all three terms on the left side of the above inequality are asymptotic to $n \log n$, 
but only two terms on the right side are;
the remaining terms on the right side are asymptotically of smaller 
order of magnitude than $n \log n$, and therefore negligible.
It turns out that the above inequality holds for $n=33$.
By the Rosser-Schoenfeld inequalities, therefore
$$p_{33}+p_{32}+p_{31}>p_{34}+p_{35}.$$
For larger $n$, $p_{n}+p_{n-1}+p_{n-2}>p_{n+1}+p_{n+2}$ must hold 
because the left side of the Rosser-Schoenfeld bounds
$n \log n + (n-1) \log(n-1) + (n-2) \log(n-2)$
grows faster than the right side
$(n+1)(\log(n+1) + \log \log(n+1)) + (n+2)(\log(n+2) + \log \log(n+2))$.
(But look out: $p_{n}+p_{n-1}+p_{n-2}-(p_{n+1}+p_{n+2})$
is neither weakly nor strictly increasing with increasing $n$,
even for $n\geq 33$.)
This $n=33$ is higher than necessary,
because it is readily verified that
$$p_{10}+p_9 + p_8=71> p_{11}+ p_{12}=68,$$
and that $p_{n}+p_{n-1}+p_{n-2}>p_{n+1}+p_{n+2}$
holds for all $n$ from 10 to 33.
Since we have sketched the proof that the inequality must continue to hold
for $n$ larger than 33,
we conclude that $p_{n}+p_{n-1}+p_{n-2}>p_{n+1}+p_{n+2}$
for all $n\ge N=10$.

\section{Special cases}

\begin{theorem}\label{thm:c1} 
For all $n=2, 3, \ldots$, we have $p_{n-1} + p_{n}\ge p_{n+1}$.
Equality holds only for $n=2$. 
\end{theorem}
We give a brief proof that is independent of Theorem \ref{thm:general}.
\begin{proof}
Loo \cite[p.\ 1880, Corollary 2.2]{Loo2011} showed that for any integer $n\ge 3$, 
there is a prime in the interval $(n,4(n+2)/3)$. 
For $n=1$ and $n=2$, 
there is a prime in the interval $(n,4(n+2)/3)$
because the corresponding intervals are $(1, 4)$ and $(2, 16/3)$.
Since $p_{n}$ is the smallest prime larger than $p_{n-1}$, 
\begin{equation*}
p_{n}\in \left(p_{n-1},\dfrac{4(p_{n-1}+2)}{3}\right). 
\end{equation*}
Likewise, 
\begin{align*}
p_{n+1}&\in \left(p_{n},\dfrac{4(p_{n}+2)}{3}\right) \\
&\subset \left(p_{n},\dfrac{4(4(p_{n-1}+2)/3+2)}{3}\right) =\left(p_{n},\dfrac{16p_{n-1}+56}{9}\right).
\end{align*}
Now 
\begin{equation*}
\dfrac{16p_{n-1}+56}{9} \le 2p_{n-1}
\end{equation*}
if and only if $28\le p_{n-1}$. So for all $p_{n-1}>28$ we have $p_{n+1} < 2p_{n-1} < p_{n-1}+p_{n}$.

To conclude the proof, it remains only verify that $p_{n-1} + p_{n}\ge p_{n+1}$
for the primes $p_{n-1}$ less than 28:
 $2+3=5, 3+5>7, 5+7>11, 7+11>13, 11+13>17, 13+17>19, 17+19>23, 19+23>29, 23+29>31$.

Because all primes $p_n$ for $n>1$ are odd and the sum of two such primes is even,
the strict inequality $p_{n-1} + p_{n}> p_{n+1}$ must hold for $n>2$. 
\end{proof}
For any positive integers $c, d$ (not necessarily distinct), 
define $N(c,d)$ to be the least finite positive integer, if it exists, 
such that $p_{n-c} + p_n \ge p_{n+d}$
for all $n \ge N(c,d)$.
Theorem \ref{thm:general} proves that $N(c,d)$ always exists.
Theorem \ref{thm:c1} proves that $N(1,1)=2$. 
Using Shevelev et al.'s result \cite{Shevelev2013} that the list of integers $k$ for which every interval 
$(kn, (k + 1)n), n > 1$, contains
a prime includes $k = 1, 2, 3, 5, 9, 14$ and no other values of $k \le 10^8$, 
we proved by extended calculations analogous to those in 
the proof of Theorem \ref{thm:c1} 
that $N(2,2)=6, N(3,3)=10, N(4,4)=11, N(5,5)=15$.

\section{Numerical results and open questions}

Define $\delta(c,d,n):=p_{n-c} + p_n - p_{n+d}$.
Then $\delta(c,d,n)$ is not always a monotonic increasing
function of $n$ even when $\delta(c,d,n)>0$.
For example, $\delta(2,2,6)=7+13-19=1$,
$\delta(2,2,7)=11+17-23=5$,
$\delta(2,2,8)=13+19-29=3$,
and $\delta(2,2,9)=17+23-31=9$.
We can also have successions of two or three (or perhaps more) 
identical values of $\delta(c,d,n)$ with given $c, d$ and increasing $n$,
e.g., $\delta(1,1,50)=p_{49}+p_{50}-p_{51}=223=\delta(1,1,51)=p_{50}+p_{51}-p_{52}$.
Is there any finite upper limit to the number of identical successive
values of $\delta(c,d,n)$ with given $c, d$ and increasing $n$?

For $c=1,\ldots,6$ and $d=1,\ldots,6$,
we calculated $\delta(c,d,n)$ numerically for all the primes less than $10^{10}$
and recorded as $M(c,d)$ (Table \ref{tab:Mcd}) 
the least $n$ such that $p_{n-c} + p_n \ge p_{n+d}$ for that $n$ and all larger observed values of $n$.
We distinguish the values $M(c,d)$ calculated numerically,
using a finite (though large) set of primes, 
from the proved values $N(c,d)$.
The first five diagonal elements $M(c,c), c=1,\ldots,5$ are consistent with the values of $N(c,c)$ proved above.
For $c < d$ in Table \ref{tab:Mcd}, usually $M(c,d)>M(d,c)$,
e.g., $M(1,3) = 6 > M(3,1)=5$, but not always: $M(5,6)=15<M(6,5)=16$. 
In Table \ref{tab:Mcd} and much larger tables not reproduced here,
for a given $c$, $M(c,d)$ is monotonically weakly increasing with $d$,
and for a given $d$, $M(c,d)$ is monotonically weakly increasing with $c$.
Is this always true?

\begin{table}[htb]
\centering
\begin{tabular}{c  r r r r r r}
\hline
     $c\downarrow, d\rightarrow$ &     1  &     2 &      3 &      4  &     5  &     6\\
\hline
\cline{1-7}
     1   &    2    &   5    &   6   &    9   &   10  &    11 \\
     2   &   3    &    6    &    8    &  10    &   10    &   12 \\
     3   &     5     &   7    &   10    &   10    &   12   &   13 \\
     4   &   5     &   9    &   10    &   11     &  13    &   14 \\
     5    &    7   &    10    &   11    &   12    &   15    &   15 \\
     6     &   7    &   11    &   12    &   14    &   16    &   16 \\
\hline
\end{tabular}
\caption{The entry in the row labeled $c$ and column headed $d$ is the numerically calculated value $M(c,d)$ 
such that, for all $n\ge M(c,d), p_{n-c} + p_n \ge p_{n+d}$ among the 
455,052,511 primes less than $10^{10}$, and such that for $n= M(c,d)-1, p_{n-c} + p_n < p_{n+d}$.
For example, $M(2,3)=8$ asserts that $p_6 + p_8 \ge p_{11}$ 
and the inequality $p_{n-2} + p_n \ge p_{n+3}$
holds for every $n\geq 8$ in this finite set of primes but not for $n=7$.
Here, as claimed, $13+19>31$ but $11+17<29$.
}
\label{tab:Mcd}
\end{table}

The first 100 values of $M(1,d), d=1,\ldots, 100$, are
2     5     6     9    10    11    12    13    14    17    17    17    20    22
24    25    25    26    26    31    31    32    32    34    35    35    38    38    41
42    44    44    47    48    48    48    49    49    52    54    55    57    62    63
63    63    64    64    64    67   67    68    68    69    69    74    74    75    76   79    81    81    81    82    84    84    87    92    93    94    94    96    98    98   99    99   100   100   100   101   102   102   102   103   104   109   112   113   115   117   117   119   120   120   122   127   128   129   129   130.
A search of the
{\it On-Line Encyclopedia of Integer Sequences} \cite[2021-04-26]{oeis} 
revealed no matching sequences.

\section{Acknowledgments}
I thank Pierre Deligne and two anonymous referees for very helpful comments, and Roseanne Benjamin for help during this work.

\end{document}